\documentclass[a4paper,11pt]{article}
\usepackage[english]{babel}
\usepackage{mathrsfs}

\usepackage{amsfonts}
\usepackage[centertags]{amsmath}
\usepackage{amsthm}
\usepackage{enumitem}
\usepackage{latexsym,amsmath,amssymb}
\usepackage{graphicx,color}
\newcommand{\ds}{\displaystyle}

\setenumerate[1]{label=(\roman*), ref=(\roman*)}
\usepackage[all]{xy}

\newtheorem{theorem}{Theorem}
\newtheorem{proposition}[theorem]{Proposition}
\newtheorem{definition}[theorem]{Definition}

\newtheoremstyle{obs}
  {3pt}
  {3pt}
  {}
  {}
  {\bfseries}
  {.}
  {.5em}
  {}
\theoremstyle{obs}
\newtheorem{remark}[theorem]{Remark}
\newtheorem{state}[theorem]{Statement}

\def\qed{\quad{$\blacksquare$}\medskip }

\def\qed{\ifvmode\removelastskip\fi
{\unskip\nobreak\hfil\penalty50\hbox{}\nobreak\hfil \hbox{\vrule
height1.2ex width1.2ex}\parfillskip=0pt \finalhyphendemerits=0
\par \smallskip}}

\title{Presymplectic high order maximum principle}
\author{\textsc{M. Barbero\textendash{}Li\~n\'an}\\
 Department of
Mathematics \& Statistics, Queen's University, \\ Jeffery Hall,
University Ave. Kingston, ON, Canada, K7L 3N6
\\[3mm]
\textsc{M. C. Mu\~noz-Lecanda}
\\
 Departamento de Matem\'atica Aplicada IV,  Universitat Polit\`ecnica de Catalunya,\\
 Edificio C-3, Campus Norte UPC,
   C/ Jordi Girona 1. 08034 Barcelona, Spain}

\begin{document}

\maketitle

\begin{abstract}
Pontryagin's Maximum Principle is an outstanding result for solving
optimal control problems by means of optimizing a specific function
on some particular variables, the so called controls. However, this
is not always enough for solving all these problems. A high order
maximum principle (Krener, 1977) must be used in order to obtain
more necessary conditions for optimality. These new conditions
determine candidates to be optimal controls for a wider range of
optimal control problems.

Here, we focus on control-affine systems. Krener's high order
perturbations are redefined following the notions introduced in
Aguilar\textendash{}Lewis (2008). A weaker version of Krener's high
order maximum principle is stated in the framework of presymplectic
geometry. As a result, the presymplectic constraint algorithm in the
sense of Gotay\textendash{}Nester\textendash{}Hinds (1979) can be
used. We establish the connections between the presymplectic
constraint algorithm and the candidates to be optimal curves
obtained from the necessary conditions in Krener's high order
maximum principle. In this paper we obtain weaker geometric
necessary conditions for optimality of abnormal solutions than the
ones in Krener (1977) and the ones in the weak high order maximum
principle. These new necessary conditions are more useful,
computationally speaking, for finding curves candidate to be
optimal. The theory is supported by describing specifically some of
the above-mentioned conditions for some mechanical control systems.


\textbf{Mathematics Subeject Classification (2000):} 49E25, 49J15, 70H05, 70H50, 93.40.
\end{abstract}

\section{Introduction}

Pontryagin's Maximum Principle (PMP) \cite{P62} is a useful result
to try to solve optimal control problems since it provides necessary
conditions for optimality. This Principle is understood as a first
order test of optimality for two reasons. If the controls take
values in the interior of the control set, the condition of
maximization of the Hamiltonian can be replaced by the fact that the
partial derivative of the Hamiltonian with respect to the controls
vanishes. On the other hand, the perturbations of the optimal curves
at a point are linearly approximated by the so-called perturbation
vectors in a neighbourhood of the point.

Despite the utility of Pontryagin's Maximum Principle it is not
always possible to determine the optimal controls from its necessary
conditions \cite{2005BulloLewis,1977Krener}. For instance, the
singular optimal curves are those curves whose controls are not
uniquely determined from the condition of maximization of the
Hamiltonian \cite{2005BulloLewis,1977Krener}. There are also the
so-called abnormal extremals that cannot be determined in general by
means of Pontryagin's Maximum Principle \cite{2005BulloLewis,LS96}.
That is why research toward a high order maximum principle started
to be developed \cite{98Bianchini,KawskiSurvey,Knobloch,1977Krener}.
Hermes \cite{Hermes} was one of the first to introduce high order
perturbations to study controllability along a fixed trajectory.
Later these perturbations were adapted to provide necessary
conditions for optimality \cite{Knobloch,1977Krener}.

Geometrically speaking, optimal control problems can be studied from
the viewpoint of symplectic geometry or of presymplectic one. The
former consists of fixing the controls in such a way that the
natural symplectic structure of the cotangent bundle of the state
space is considered \cite{2004Agrachev,2008PMPMiguelMaria,S98Free}.
The latter consists of using the presymplectic structure
\cite{2007MiguelMaria} of the above cotangent bundle times the
control set. In the presymplectic formalism, we are able to state a
weaker version of Krener's high order maximum principle in the same
way as stated the weaker version of Pontryagin's Maximum Principle
in \cite{2007MiguelMaria}. Some of the necessary conditions for
optimality coming from the weak high order maximum principle give a
geometric meaning to previous techniques used in the literature
\cite{KawskiSurvey,Knobloch,1977Krener}.

The paper is organized as follows. In Section \ref{SBackground} we
review the notion of an optimal control problem using two equivalent
statements. Section \ref{HOMP} contains all the necessary elements
to state Krener's high order Maximum Principle \cite{1977Krener} for
control-affine systems by using the constructions in
\cite{MTNSCesar}. In Section \ref{Sweak} a weak high order maximum
principle is stated. The main result of this paper consists of
establishing a connection between the above-mentioned weak principle
and the classical high order maximum principle \cite{1977Krener} as
described in Sections \ref{Discuss} and \ref{ExHOMP} under the
assumption of having an abnormal optimal solution. In order to
clarify ideas, we conclude with a detailed description of some
elements that appear in Section \ref{Discuss} for some mechanical
control systems in Section \ref{ExHOMP}.

Before concluding this introduction we want to point out that most
of the theory in the first part of Section \ref{HOMP} and Section
\ref{Sweak} can be rewritten for any control system. We restrict
ourselves to control-affine systems because the aim of this paper is
to establish a relationship between the above-mentioned
presymplectic constraint algorithm and the high order maximum
principle as explained in Section \ref{Discuss} for abnormal optimal
solutions. This relationship is obtained by identifying brackets of
vector fields in the control system with high order perturbation
vectors. To our best knowledge, it is not known how to deal with
that when the control systems are not control-affine.

\section{Background in optimal control theory}\label{SBackground}

First, let us define geometrically the notion of control system.

\begin{definition}\label{DefControlSyst} Let $M$ be an $m$-dimensional manifold and $U$ be a
subset of $\mathbb{R}^k$. A \textbf{control system} on $M$ is a
smooth vector field $X$ defined along the projection $\pi\colon M
\times U \rightarrow M$. A \textbf{trajectory} or an
\textbf{integral curve of the control system $X$} is a curve
$(\gamma,u)\colon I \subset \mathbb{R} \rightarrow M\times U$ such
that $\gamma$ is absolutely continuous, $u$ is measurable and
bounded, and $\dot{\gamma}(t)=X(\gamma(t),u(t))$ a.e. for $t\in I$.
\end{definition}

The set of control systems is denoted by $\mathfrak{X}(\pi)$; that
is, the set of smooth vector fields defined along the projection
$\pi$. The curve $u\colon I \rightarrow U$ is called the
\textit{control}.

Given a control system $X\in \mathfrak{X}(\pi)$, we are usually
interested in the set of points that can be reached from an initial
point through trajectories $(\gamma,u)\colon I\subset \mathbb{R}
\rightarrow M\times U$ of $X$, where $I=[a,b]$. In this regard, the
following definition is essential.
\begin{definition}\label{ReachSet} Let $M$ be a manifold, $U$ be a set in
$\mathbb{R}^k$ and $X$ be a smooth vector field along the projection
$\pi\colon M\times U \rightarrow M$. The \textbf{reachable set
 from $x_0\in M$ at time $T\in I$} is the set
of points described by \begin{eqnarray*} {\mathcal
R}(x_0,T)=\left\{x\in M  \right. & | & \textrm{ there exists }
(\gamma,u)\colon [a,b] \rightarrow M\times U \textrm{ such that } \\
& & \left. \dot{\gamma}(t)=X(\gamma(t),u(t)) \textrm{ a.e.}, \;
\gamma(a)=x_0, \; \gamma(T)=x \right\}.\end{eqnarray*} The
\textbf{reachable set from a point $x_0\in M$ up to time $T$} is
\[{\mathcal R}(x_0,\leq T)=\bigcup_{a\leq t \leq T} {\mathcal
R}(x_0,t).\]
\end{definition}

The topological properties of the reachable set are closely related
to different properties of the control system such as reachability,
controllability, see \cite{2005BulloAndrewBook} for more details.
Moreover, in the sequel we briefly describe why the linear
approximation of the reachable sets is the key point to prove
Pontryagin's Maximum Principle \cite{P62} and also Krener's high
order maximum principle \cite{1977Krener}.

\subsection{Optimal control problem}

Let us define an optimal control problem associated to the control
system in Definition \ref{DefControlSyst}. Given ${\mathcal F}
\colon M\times U \rightarrow \mathbb{R}$, consider the functional
$${\cal S}[\gamma,u]=\int_I {\mathcal F}(\gamma(t),u(t))\, {\rm d}t$$
defined on curves $(\gamma, u)$ with a compact interval as domain.
The function ${\mathcal F} \colon M \times U \rightarrow \mathbb{R}$
is continuous on $M \times U$ and continuously diffe\-rentia\-ble
with respect to $M$ on $M \times U$.

If we fix the control $u$, the vector field $X\in \mathfrak{X}(\pi)$
can be rewritten as a time-dependent vector field $X^{\{u\}}\colon
I\times M \rightarrow TM$, $X^{\{u\}}(t,x)=X(x,u(t))$.

\begin{state}\label{OCP} \textbf{(Optimal Control Problem, OCP)}
Given the elements $M$, $U$, $X$, ${\mathcal F}$, $I=[a,b]$ and the
endpoint conditions $x_a$, $x_b\in M$, consider the following
problem.

Find $(\gamma,u)$ such that
\begin{itemize}
\item[(1)] $\gamma(a)=x_a$, $\gamma(b)=x_b$ (endpoint conditions),
\item[(2)] $\gamma$ is an integral curve of $X^{\{u\}}$, i. e.
$\dot{\gamma}(t)=X^{\{u\}}(t,\gamma(t))$, for a.e. $t\in I$, and
\item[(3)] ${\cal S}[\gamma,u]$ is minimum over all curves
satisfying  (1) and (2) (minimal condition).
\end{itemize}
\end{state}
The function ${\mathcal F}$ is called the \textit{cost function} of
the problem. A \textit{solution to the OCP} is a curve $(\gamma,u)$
satisfying conditions (1)-(3) in Problem \ref{OCP} such that
$\gamma$ is absolutely continuous and $u$ is measurable and bounded.
These problems are also known in the literature as Lagrange
problems.

\subsection{Extended optimal control problem}
A usual technique in optimal control theory is to consider an
equivalent problem to Problem \ref{OCP}, but defined on the extended
manifold $\widehat{M}=\mathbb{R} \times M$. Let $\widehat{X}$ be the
following vector field defined along the projection
$\widehat{\pi}\colon \widehat{M} \times U \rightarrow \widehat{M}$,
$$\widehat{X}(x^0,x,u)={\mathcal F}(x,u) {\partial}/{\partial x^0}|_{(x^0,x,u)}
+ X(x,u),$$ where $x^0$ is the natural coordinate on $\mathbb{R}$.
In this new manifold it is easier to identify the direction of
decreasing of the cost function. Now the integral curves of the
control system $\widehat{X}$ are curves $(\widehat{\gamma},u)=((x^0
\circ \widehat{\gamma}, \gamma), u) \colon I \rightarrow \widehat{M}
\times U$ such that $\widehat{\gamma}$ is absolutely continuous and
$u$ is measurable and bounded.

\begin{state}\label{stateEOCP}\textbf{(Extended Optimal Control Problem, $\widehat{OCP}$)}
Given the elements in Problem \ref{OCP}, $\widehat{M}$ and
$\widehat{X}$, consider the fo\-llo\-wing problem.

Find $(\widehat{\gamma},u)=(\gamma^0,\gamma,u)$ such that
\begin{itemize}
\item[(1)] $\widehat{\gamma}(a)=(0,x_a)$, $\gamma(b)=x_b$ (endpoint conditions),
\item[(2)] $\widehat{\gamma}$ is an integral curve of $\widehat{X}^{\{u\}}$, i.e.  $\dot{\widehat{\gamma}}(t)=\widehat{X}^{\{u\}}(t,\widehat{\gamma}(t))$,
for a.e. $t\in I$, and
\item[(3)] $\gamma^0(b)$ is minimum over all curves
satisfying  (1) and (2) (minimal condition).
\end{itemize}
\end{state}

A \textit{solution to the $\widehat{OCP}$} is a curve
$(\widehat{\gamma},u)$ satisfying conditions (1)-(3) in Problem
\ref{stateEOCP} such that $\widehat{\gamma}$ is absolutely
continuous and $u$ is measurable and bounded. In the literature this
problem is also known as Mayer problem.

\section{High order Maximum Principle for control-affine
systems} \label{HOMP}

First order necessary conditions for optimality do not always
provide enough information for finding the optimal solution. High
order variations of a curve depending on parameters must be studied,
as is done for instance in
\cite{98Bianchini,BS93,72GabKir,KawskiSurvey,Knobloch,1977Krener},
in order to obtain high order necessary conditions for optimality.

Given a fixed reference trajectory $(\gamma,u)$ defined in $I\subset
\mathbb{R}$, we use the description of high order variations for
control-affine systems considered in \cite{MTNSCesar} to study
controllability at a point $x_0\in M$  in order to obtain the
variations along a reference trajectory defined in
\cite{1977Krener}.

Let $\Sigma=(M,\{X_0,X_1,\dots,X_k\},U)$ be a {\em control-affine
system} whose trajectories $(\gamma,u)$ on $M\times U$ satisfy
\[\dot{\gamma}(t)=X_0(\gamma(t))+\sum_{c=1}^ku^cX_c(\gamma(t)).\]
The vector field $X_0$ is usually called the {\em drift vector
field} and $X_1,\dots,X_k$ are the {\em control or input vector
fields}. For the sake of simplicity all the vector fields are
assumed to be complete from now on.
\begin{definition}\label{highpert} Let $\mathbf{\xi}=(\xi_1,\ldots,\xi_r)$ be a finite sequence of vector
fields on $M$ where $\xi_i(x)\in \{X_0(x)+\sum_{c=1}^ku^cX_c(x) \mid
u\in U\}$ for each $i\in \{1,\ldots, r\}$ and $x\in M$.
\begin{enumerate}
\item The flow associated with $\mathbf{\xi}$ called $\mathbf{\xi}$\textbf{-flow} is the map
$\Phi^{\mathbf{\xi}}\colon \mathbb{R}^r\times M \rightarrow M$ given
by
\begin{equation}\label{gammapert}
\Phi^{\mathbf{\xi}}(\mathbf{t},x)=\left(\Phi_{t_r}^{\xi_r}\circ
\dots \circ \Phi_{t_1}^{\xi_1}\right)(x),\end{equation}
 where
$\mathbf{t}=(t_1,\ldots,t_r)$. For a fixed $x\in M$ define
$\Phi^{\mathbf{\xi}}_x\colon \mathbb{R}^r\rightarrow M$ by
$\Phi^{\mathbf{\xi}}_x(\mathbf{t})=\Phi^{\mathbf{\xi}}(\mathbf{t},x)$,
and for a fixed $\mathbf{t}\in \mathbb{R}^r$ define
$\Phi^{\mathbf{\xi}}_{\mathbf{t}}\colon M\rightarrow M$ by
$\Phi^{\mathbf{\xi}}_{\mathbf{t}}(x)=\Phi^{\mathbf{\xi}}(\mathbf{t},x)$.
\item A \textbf{positive $r$-end-time
variation} is a smooth map $\mathbf{\tau}\colon
[0,\infty)\rightarrow \mathbb{R}^r_{\geq 0}$ with the property
$\tau(0)=\mathbf{0}$. The set of all the positive $r$-end-time
variation is denoted by ${\rm ET}^+_r$.
\end{enumerate}
\end{definition}

As mentioned, the variations in Definition \ref{highpert} are useful
for studying particular notions of controllability as shown in
\cite{MTNSCesar}. However, to study optimality the variations
described in \cite{1977Krener} are necessary. We describe them
similarly to the ones in Definition \ref{highpert} using
some usual notation from the theory of jet bundles \cite{Sa-89}. 

\begin{definition}\label{Def-Krener-variations}
Let $(\gamma,u)$ be a reference trajectory associated with the
vector field $\xi_0$ and $\mathbf{\xi}=(\xi_1,\ldots,\xi_r)$ be a
finite sequence of vector fields on $M$ where $\xi_i(x)\in
\{X_0(x)+\sum_{c=1}^ku^cX_c(x) \mid u\in U\}$ for each $i\in
\{0,\ldots, r\}$, $x\in M$, $\tau \in {\rm ET}^+_r$.
\begin{enumerate}
\item A \textbf{$(r+2)$-end-time variation}
is a smooth map  $\tau_2:[0,\infty)\rightarrow \mathbb{R}^2\times
\mathbb{R}^r_{\geq 0}$ such that $\tau_2(s)=(q_1(s), q_2(s),
\tau(s))$ and $\tau_2(0)=\mathbf{0}$.

\item The $(\xi,\tau_2)$-\textbf{variation along $\gamma$ at time $t_0$} is the
curve $\nu_{\xi,\tau_2}\colon [0,\infty) \rightarrow M$ given by
\begin{equation}\label{Eq-Var-2}\nu_{\xi,\tau_2}(s)=(\Phi^{\xi_0}_{q_2(s)}\circ
\Phi^\xi_{\tau(s)}\circ
\Phi^{\xi_0}_{q_1(s)})(\gamma(t_0)).\end{equation}
\item
If $\tau_2$ is $(r+2)$-end-time variation, the \textbf{order} of the
pair $(\xi,\tau_2)$ at $\gamma(t_0)\in M$, denoted ${\rm
ord}_{\gamma(t_0)}(\xi,\tau_2)$, is the smallest positive integer
$l$ such that there exists $\epsilon>0$ satisfying
\begin{equation}\label{Eq-order-general}
j^l(\nu_{\xi,\tau_2})(s)\neq 0_{\nu_{\xi,\tau_2}(s)}, \quad
|s|<\epsilon.\end{equation} If no such $s$ exists then the order is
set to $\infty$. \item  If ${\rm ord}_{\gamma(t_0)}(\xi,\tau_2)=l$,
then the \textbf{$(\xi,\tau_2)$-infinitesimal variation at
$\gamma(t_0)$} is the tangent vector
\[V_{\xi,\tau_2}=j^l(\nu_{\xi,\tau_2})(0).\]
A $(\xi,\tau_2)$-infinitesimal variation of order $l$ is also called
a \textbf{high order elementary perturbation vector at $\gamma(t_0)$
of order $l$}.
\end{enumerate}\end{definition}

Note that the notion of order is defined by a condition that must be
satisfied in a neighborhood of $0$, see (\ref{Eq-order-general}). As
the controls are assumed to be measurable, this condition is
necessary to guarantee that perturbations at different times can be
summed as explained in \cite{1977Krener}.

These infinitesimal variations include the elementary perturbation
vectors in the classical Pontryagin's Maximum Principle \cite{P62}.
They are defined as follows. Take $\xi=(\xi_1)$ such that
$\xi_1(x)=X_0(x)+\sum_{c=1}^ku^c_1X_c(x)$ being $u_1\in U$,
$\tau(s)=l_1s$ for $l_1\in \mathbb{R}^+$, $q_1(s)=-l_1s$,
$q_2(s)=0$. Then the $(\xi,\tau_2)$-variation along $\gamma$
corresponding with $\xi_0(x)=X_0(x)+\sum_{c=1}^ku^cX_c(x)$ at time
$t_0$ is
\[\nu_{\xi,\tau_2}(s)=(\Phi^{\xi_1}_{l_1s}\circ \Phi^{\xi_0}_{-l_1
s})( \gamma(t_0)).\] The $(\xi, \tau_2)$-infinitesimal variation is
given by
\begin{align}
\ds{\left.\frac{{\rm d}}{{\rm d}s}\right|_{s=0}\nu_{\xi,\tau_2}(s)}
&\ds{ = \left.\frac{{\rm d}}{{\rm
d}s}\right|_{s=0}\Phi^{\xi_1}(\tau(s),\Phi^{\xi_0}(q_1(s),\gamma(t_0)))}
\nonumber
\\& \ds{+\left.\frac{\partial
\Phi^{\xi_1}}{\partial
x^i}\right|_{s=0}(\tau(s),\Phi^{\xi_0}(q_1(s),\gamma(t_0)))\left.\frac{{\rm
d}}{{\rm d}
s}\right|_{s=0}\left(\Phi^{\xi_0}\right)^i(q_1(s),\gamma(t_0))}
\nonumber
\\ &\ds{= l_1\xi_1(\gamma(t_0))+\left.\frac{\partial \Phi^{\xi_1}}{\partial
x^i}\right|_{s=0}(0,\Phi^{\xi_0}(q_1(s),\gamma(t_0)))\xi_0^i(\gamma(t_0))(-l_1)}
\nonumber \\&\ds{=l_1(\xi_1(\gamma(t_0))-\xi_0(\gamma(t_0)))
=\sum_{c=1}^kl_1(u^c_1-u^c(t_0))X_c(\gamma(t_0)),}
\label{Eq-Krener-PMP} \end{align} because $\Phi^{\xi_1}_0$ is the
identity map. The last equality in (\ref{Eq-Krener-PMP}) is exactly
Pontryagin's elementary perturbation vectors for control-affine
systems.

\begin{remark}\label{Rgeneral} All the theory developed in this section so far is true for
general control systems apart from the last equality in
(\ref{Eq-Krener-PMP}), only true for control-affine systems.
\end{remark}

\begin{definition}\label{Def-High-Pert-Vect}
Let $(\gamma,u)$ be a reference trajectory associated with $\xi_0$,
$t_0\in I$. The \textbf{set of all the high order elementary
perturbation vectors at $\gamma(t_0)$ of order $l$} is
\[\begin{array}{rcl} \mathscr{V}^l_{\Sigma_2}(\gamma(t_0)) & =& \left\{V_{\xi,\tau_2}\in T_{\gamma(t_0)}M \mid
\tau_2\in {\mathcal C}^{\infty}(\mathbb{R},\mathbb{R}^2)\times
ET^+_r, \, {\rm ord}_{\gamma(t_0)}(\xi,\tau_2)=l, \right. \\ &&
\left. \mathbf{\xi}=(\xi_1,\ldots,\xi_r), \; \xi_i\in
\{X_0+\sum_{c=1}^ku^cX_c \mid u\in U\}\right. \\ && \left.\mbox{ for
each } i=1,\dots,r,\; r\in \mathbb{Z}_{\geq 1}\right\}.\end{array}\]
The \textbf{set of all the high order elementary perturbation
vectors at $\gamma(t_0)$} is
\[\mathscr{V}_{\Sigma_2}(\gamma(t_0))=\bigcup_{l\geq 1}
\mathscr{V}^l_{\Sigma_2}(\gamma(t_0)).\]
\end{definition}

\begin{remark} \label{remark} If the control system is control-affine and $0\in {\rm
int} \, U$ then (\ref{Eq-Krener-PMP}) guarantees that \[{\rm
cone}({\rm conv} \{X_1(\gamma(t)),\dots,X_k(\gamma(t))\}) \subseteq
\mathscr{V}^1_{\Sigma_2}(\gamma(t)),\] where ${\rm conv}$ denotes
convex hull.
\end{remark}

\begin{proposition} The set $\mathscr{V}_{\Sigma_2}(\gamma(t))$ is a convex
cone.
\end{proposition}

\proof The set $\mathscr{V}_{\Sigma_2}(\gamma(t))$ is a convex cone
if it is closed under addition and closed under
$\mathbb{R}_{>0}$-multiplication.

Then \cite[Lemma 3.1]{1977Krener} guarantees that the order of any
variation can be shifted upward and \cite[Lemma 3.4]{1977Krener}
assures that the sum of variations of different order is a variation
of higher order. Thus $\mathscr{V}_{\Sigma_2}(\gamma(t))$ is closed
under addition. To prove that $\mathscr{V}_{\Sigma_2}(\gamma(t))$ is
closed under $\mathbb{R}_{>0}$-multiplication, the time must be
re-scaled by using a different $(r+2)$-end-time variation. \qed

\begin{definition}\label{Def-High-Cone} Let $(\gamma,u)$ be a reference trajectory associated with
$\xi_0$, $t\in I$. The \textbf{high order tangent perturbation cone
$K_t$ at $\gamma(t)$} is the smallest closed convex cone in
$T_{\gamma(t)}M$ that contains all the displacements by the flow of
$\xi_0$ of all the $(\xi,\tau_2)$-infinitesimal variations along
$\gamma$ at all Lebesgue times smaller than $t$, i.e.,
\begin{equation}\label{Eq-Kt}K_t=\overline{{\rm conv}\left(\bigcup_{\substack{ a<t_0\leq t \\
t_0 \textrm{ is a Lebesgue time}}} (\Phi_{(t-t_0)}^{\xi_0})_*
(\mathscr{V}_{\Sigma_2}(\gamma(t_0)))\right)}.\end{equation}
\end{definition}
\begin{remark} The notion of Lebesgue time comes from considering
Carath\'eodory vector fields, see for instance \cite{55Coddington}
for more details.
\end{remark}

 This cone generalizes
the notion of tangent perturbation cone introduced in Pontryagin's
Maximum Principle \cite{2004Agrachev,2008PMPMiguelMaria,P62} so as
to state and prove the high order maximum principle, see Theorem
\ref{HighPMP}.

The high order tangent pertubation cone may be understood as an
approximation of the reachable set along a reference trajectory that
is perturbed as shown in the following result.

\begin{proposition} \label{lemma2} Let $(\gamma,u)$ be a reference trajectory and $t \in (a,b]$. If $v$ is a nonzero vector
in the interior of the high order tangent perturbation cone $K_t$,
then there exist $\epsilon>0$ and a
$(\xi,\tau_2)$-variation of order $l$ 
such that
\[ \nu_{\mathbf{\xi},\mathbf{\tau}_2}(s)
=\gamma(t)+v\frac{s^l}{l!}+o(s^l)\quad {\rm for} \;\;
0<s<\epsilon.\]
\end{proposition}

\proof It is part of the proof of Theorem 3.6 in \cite[page 266
]{1977Krener}. \qed

We have constructed the above perturbations on the manifold $M$ for
the sake of simplicity, but in order to state the high order maximum
principle these perturbations must be considered for the extended
control system defined on $\widehat{M}=\mathbb{R}\times M$. On
$\widehat{M}$ the perturbations are defined analogously as the ones
in Definition \ref{Def-Krener-variations}. Note that in general a
control-affine system on $M$ is not control-affine on $\widehat{M}$.
However, as mentioned in Remark \ref{Rgeneral}, all the above
constructions and results are true for any control system apart from
Remark \ref{remark}.

To solve optimal control problems is difficult. That is why
Pontryagin's Maximum Principle \cite{P62} transforms them into
Hamiltonian problems that provide more conditions to characterize
optimality. To state Krener's maximum principle \cite{1977Krener} we
need to define a Hamiltonian system associated with the control
system $\widehat{X}$. First, consider the cotangent bundle $T^{\ast}
\widehat{M}$ with its natural symplectic structure that will be
denoted by $\Omega$. If
$(\widehat{x},\widehat{p})=(x^0,x,p_0,p)=(x^0, x^1, \ldots ,
x^m,p_0, p_1, \ldots , p_m)$ are local natural coordinates on $T^*
\widehat{M}$, the form $\Omega$ has as local expression $\Omega={\rm
d}x^0\wedge {\rm d}p_0+{\rm d}x^i\wedge {\rm d}p_i$. For each $u \in
U$, $H^u \colon T^* \widehat{M} \rightarrow \mathbb{R}$ is the
Hamiltonian function defined by
$$H^u(\widehat{p})=H(\widehat{p},u)=
\langle \widehat{p},\widehat{X}(\widehat{x},u) \rangle =p_0
{\mathcal F}(x,u) +\sum_{i=1}^m p_i f^i(x,u),$$ where
$\widehat{p}\in T^*_{\widehat{x}}\widehat{M}$ and
$\widehat{X}={\mathcal F}\partial /
\partial x^0+f^i\partial /
\partial x^i$. The tuple $(T^*\widehat{M}, \Omega, H^u)$ is a Hamiltonian
system. The associated Hamiltonian vector field
$\widehat{X}_H^{\{u\}}$ on $T^*\widehat{M}$ defined by
$\widehat{X}_H^{\{u\}}(\widehat{p})=\widehat{X}_H(\widehat{p},u)$
satisfies the equation
\begin{equation*} i_{\widehat{X}_H^{\{u\}}} \Omega ={\rm
d}H^u.\end{equation*} It should be noted that
$\widehat{X}_H^{\{u\}}=\left(\widehat{X}^{\{u\}}\right)^{T^*}$; that
is, $\widehat{X}_H^{\{u\}}$ is the cotangent lift of
$\widehat{X}^{\{u\}}$. Hence we get a family of Hamiltonian systems
parameterized
 by $u$ and given by
$H^u$.

With the above constructions, the Hamiltonian system in mind and the
following definition, we are able to state the high order maximum
principle \cite{1977Krener}.

\begin{definition} Let $C$ be a cone with
vertex at $0\in T_xM$. A \textbf{supporting hyperplane to $C$ at
$0$} is a hyperplane such that $C$ is contained in one of the
half-spaces defined by the hyperplane.
\end{definition}

\begin{theorem}\label{HighPMP} \textbf{(High order maximum principle, HOMP)}
Let $(\widehat{\gamma},u)\colon I \rightarrow \widehat{M}\times U$
be a solution to $\widehat{OCP}$, Problem \ref{stateEOCP}, and let
$H$ be the Hamiltonian function defined by
\begin{equation}\label{H}\begin{array}{rcl} H\colon T^*\widehat{M} \times U &\longrightarrow & \mathbb{R} \\
(\widehat{p},u)& \longmapsto & H(\widehat{p},u)=\langle \widehat{p},
\widehat{X}(\widehat{x},u) \rangle\end{array}\end{equation} where
$\widehat{p}\in T_{\widehat{x}}^*\widehat{M}$. Then there exists
$\widehat{\lambda}\colon I \rightarrow T^*\widehat{M}$ along
$\widehat{\gamma}$ such that:
\begin{enumerate}
\item $\widehat{\lambda}$ is an integral curve of a
Hamiltonian vector field $\widehat{X}_H^{\{u\}}$ on
$T^*\widehat{M}$, i.e. solution to Hamilton's equations
\begin{equation}\label{Hamiltoneq}
i_{\widehat{X}_H^{\{u\}}}\Omega={\rm d}H^{\{u\}};
\end{equation}
\item \begin{description}
\item[(a)] $H\left(\widehat{\lambda}(t),u(t)\right)= \sup_{w \in U} H(\widehat{\lambda}(t),
w)$ almost everywhere;
\item[(b)] $\sup_{w \in U} H(\widehat{\lambda}(t),
w)$ is constant everywhere;
\item[(c)] for every $t\in I$, $\ker \widehat{\lambda}(t)$ is a supporting hyperplane to $\widehat{K}_t$ at $0$;
\item[(d)] $\widehat{\lambda}(t)\neq 0 \in
T^*_{\widehat{\gamma}(t)} \widehat{M}$ for each $t\in [a,b]$;
\item[(e)] $\lambda_0$ is constant and $\lambda_0 \leq 0$.
\end{description}
\end{enumerate}
\end{theorem}

\proof It follows similarly to the proof of Pontryagin's Maximum
Principle \cite{2004Agrachev,2008PMPMiguelMaria,P62} having in mind
Proposition \ref{lemma2}. For more details see proof in
\cite[Section 3]{1977Krener}. \qed

The condition that does not appear in the classical Pontryagin's
Maximum Principle (PMP) is item (2c). The existence of this
supporting hyperplane is not explicitly stated in PMP, though it is
also fulfilled. Thus all the constructions related with the cones
are not necessary to state the classical PMP, but they are
completely necessary to prove it. After Definition
\ref{Def-Krener-variations} we have already proved how to relate
$(\xi,\tau_2)$-variations of order $1$ with the elementary
perturbation vectors used in the classical PMP, see
(\ref{Eq-Krener-PMP}). For these vectors the existence of a
supporting hyperplane is equivalent to
\[\begin{array}{ll} {}&\langle \widehat{p},\widehat{X}(\widehat{\gamma}(t),w(t))-\widehat{X}
(\widehat{\gamma}(t),u(t))\rangle \leq 0 \quad \forall \; w \;
\textrm{s.t. }  w(t)\in U,\; \; \forall \; t\in I,  \\  \mbox{ that
is, } & \langle
\widehat{p},\widehat{X}(\widehat{\gamma}(t),w(t))\rangle \leq
\langle \widehat{p},\widehat{X} (\widehat{\gamma}(t),u(t))\rangle
\quad  \forall \; w \; \textrm{s.t.  }  w(t)\in U ,\; \forall \;
t\in I.\end{array}\] In other words, the existence of a supporting
hyperplane to the first order tangent perturbation cone is
equivalent to the classical condition of maximization in PMP
corresponding with (2a) in Theorem \ref{HighPMP}. That is why in PMP
the supporting hyperplane is not mentioned in the statement.

On the contrary, in the high order version of this Principle more
perturbations are considered. The usual necessary conditions in PMP
appear in Theorem \ref{HighPMP}, but also some extra conditions as
given in (2c). These extra conditions consist of determining a
hyperplane that
 supports a greater cone than the first order tangent perturbation cone. Remember that the need of a high order maximum
principle comes from the need to obtain more conditions to determine
the optimal controls.

Observe that HOMP guarantees the existence of a covector along the
optimal curve, but it does not say anything about the uniqueness of
the covector, see Definition \ref{definextremal} for the different
kind of curves. Indeed, the key point in the proof of HOMP
\cite{1977Krener} consists of choosing an initial condition for
$\widehat{\lambda}$ to integrate (\ref{Hamiltoneq}) in such a way
that the resulting momentum satisfies all the necessary conditions
in Theorem \ref{HighPMP}. Moreover, as $\widehat{\gamma}$ is an
optimal trajectory, the momentum must separate the cone
$\widehat{K}_t$ and the direction of decreasing of the cost
function. Note that not all the supporting hyperplanes fulfill
necessarily this separating condition. The non existence of a
supporting hyperplane being also separating would immediately
contradict the hypothesis of optimality, see proof in
\cite{2004Agrachev,2008PMPMiguelMaria,1977Krener,P62} for more
details.

\begin{remark} We always refer to the maximum instead of supremum because the set of
necessary conditions in HOMP guarantee that the supremum is
achieved, so it can be called maximum.
\end{remark}

\begin{definition}  \label{definextremal}
A curve $(\widehat{\gamma}, u)\colon [a,b] \rightarrow \widehat{M}
\times U$ for $\widehat{OCP}$, Problem \ref{stateEOCP}, is
\begin{enumerate}
\item an \textbf{extremal} if there exists $\widehat{\lambda}\colon [a,b] \rightarrow
T^* \widehat{M}$ along $\widehat{\gamma}=\pi_{\widehat{M}} \circ
\widehat{\lambda}$ being $\pi_{\widehat{M}}\colon
T^*\widehat{M}\rightarrow \widehat{M}$ and $(\widehat{\lambda},
u)\colon [a,b]\rightarrow T^* \widehat{M}\times U$ satisfies the
necessary conditions of HOMP. The pair $(\widehat{\lambda}, u)$ is
called \textbf{biextremal for $\widehat{OCP}$};
\item a \textbf{normal extremal}
if it is an extremal with $\lambda_0=-1$;
\item an \textbf{abnormal extremal} if it is an extremal with
$\lambda_0=0$;
\item a \textbf{strictly abnormal
extremal} if it is not a normal extremal, but it is abnormal.
\end{enumerate}
\end{definition}

\subsection{Weak high order maximum principle}\label{Sweak}

From now on the control set $U$ is assumed to be an open set. In the
sequel, we use a presymplectic framework to state a weaker high
order Maximum Principle, as was already considered for PMP in
\cite{2007MiguelMaria}. We will show that the presymplectic
framework provides some useful techniques to characterize optimality
in a geometric way. Some of the new conditions obtained are easier
to deal with than to compute the $(\xi,\tau_2)$-infinitesimal
variations in order to find a supporting hyperplane to all of them.

Observe that $(T^*\widehat{M}\times U, \Omega, H)$ is a
presymplectic Hamiltonian system with $\Omega$ being the pullback of
the canonical 2\textendash form on $T^*\widehat{M}$ through
$\pi_1\colon T^*\widehat{M}\times U \rightarrow T^*\widehat{M}$ and
$H$ the function in (\ref{H}).

\begin{theorem}\label{weakHPMP} \textbf{(Weak high order maximum principle, WHOMP)}
Let $(\widehat{\gamma},u)\colon I \rightarrow \widehat{M}\times U$
be a solution to $\widehat{OCP}$, Problem \ref{stateEOCP}, and let
$H$ be the Hamiltonian function defined in (\ref{H}). Then there
exists $\widehat{\lambda}\colon I \rightarrow T^*\widehat{M}$ along
$\widehat{\gamma}$ such that:
\begin{enumerate}
\item $(\widehat{\lambda},u)$ is an integral curve of a
Hamiltonian vector field $\widehat{X}_H$ such that
\begin{equation}\label{preeq}
i_{\widehat{X}_H}\Omega={\rm d}H, \; {\rm i.e.}\;  \;
i_{(\dot{\widehat{\lambda}}(t),\dot{u}(t))}\Omega={\rm
d}{H}(\widehat{\lambda}(t),u(t));
\end{equation}
\item \begin{description}
\item[(a)] $H(\widehat{\lambda}(t),
u(t))$ is constant almost everywhere;
\item[(b)] for
every $t\in I$, $\ker \widehat{\lambda}(t)$ is a supporting
hyperplane to $\widehat{K}_t$ at $0$;
\item[(c)] $\widehat{\lambda}(t)\neq 0 \in
T^*_{\widehat{\gamma}(t)} \widehat{M}$ for each $t\in [a,b]$;
\item[(d)] $\lambda_0$ is constant and $\lambda_0 \leq 0$.
\end{description}
\end{enumerate}
\end{theorem}

As $\Omega$ is degenerate, (\ref{preeq}) does not necessarily have
solution on the entire manifold $T^*\widehat{M}\times U$. It may
have a solution, see \cite{1978Gotay}, if we restrict the equation
to the {\sl primary constraint submanifold} defined by
\[N_0=\left\{(\widehat{\lambda},u)\in T^*\widehat{M}\times U \mid i_v \, {\rm d}{H}=0, \quad {\rm for }\; v\in \ker
\Omega_{(\widehat{\lambda},u)} \right\}.\] Locally,
$\ds{N_0=\left\{(\widehat{\lambda},u)\in T^*\widehat{M}\times U \mid
\ds{\frac{\partial H}{\partial u^c}(\widehat{\lambda},u)=0} , \quad
c=1,\ldots, k\right\}}$. If $X_0$ is a solution to the presymplectic
equation, then $X^{N_0}=X_0+\ker \, \Omega$. At this point, a
presymplectic constraint algorithm in the sense given in
\cite{1978Gotay} starts. The stabilization process gives inductively
$N_1=\{(\widehat{\lambda},u)\in N_0 \, | \, \exists \, X\in X^{N_0}
\, , \, X(\widehat{\lambda},u)\in T_{(\widehat{\lambda},u)} N_0\}$,
$N_2,\ldots N_r$ such that $N_{r}\subseteq N_{r-1}$ and finishes
when $N_i=N_{i-1}$ for some $i\in \mathbb{N}$. The algorithm just
imposes tangency conditions in order to find a submanifold where the
optimal trajectories live and stay. This algorithm  has been adapted
when the sets are not submanifolds \cite{Eduardo}, thus this method
is not that restrictive.

Observe that $N_0$ is defined implicitly by a necessary condition
for the Hamiltonian to have an extremum over the controls as long as
$U$ is an open set. In Theorem \ref{HighPMP} the Hamiltonian is
equal to the maximum of the Hamiltonian over the controls.
Therefore, Theorem \ref{weakHPMP} is just a weaker version of
Theorem \ref{HighPMP} and the proof is straightforward from Theorem
\ref{HighPMP} and the explanation given about the presymplectic
constraint algorithm.

For normal extremals, $p_0=-1$, and for cost functions quadratic on
the controls, the presymplectic algorithm stops at $N_0$. There the
controls can already be determined. More discussion about
 how the algorithm behaves for abnormality is developed in the following section.

\subsection{Tangent perturbation cone versus
presymplectic constraint algorithm}\label{Discuss}

The high order tangent perturbation cone in Definition
\ref{Def-High-Cone} gives directions that approximate the
perturbations of the optimal trajectory in the sense given in
Proposition \ref{lemma2}.
On the other hand, the submanifolds that
appear in the presymplectic constraint algorithm live in
$T^*\widehat{M} \times U$ and they are the zero sets of a family of
vector fields $Z_j\colon \widehat{M}\times U \rightarrow
T\widehat{M}$, $j\in \mathbb{N}$ as is proved in the sequel.

Our goal now is to establish a relationship between $\widehat{K}_t$
and the sequence of constraint submanifolds $N_i$. From now on, we
focus on abnormal optimal solutions because in general the normal
optimal solutions admit a lift to the cotangent bundle
$T^*\widehat{M}$ that lies in the primary constraint submanifold and
so the algorithm stops in the first step. For abnormal extremals,
the existence of a supporting hyperplane $\ker \widehat{\lambda}$ to
$\widehat{K}_t$ is equivalent to the existence of a supporting
hyperplane given by $\ker \lambda$ to $K_t$ without considering the
extended manifold.

\begin{theorem}\label{Prop-Relation-affine} Let $U$ be an open set in $\mathbb{R}^k$,
 $\Sigma=(M,\{X_0,\ldots,X_k\},U)$ be a control-affine system
such that $0\in {\rm int}\, U$  and $(\gamma,u)$ be an abnormal
optimal solution associated with $\xi_0$. Then there exists $\lambda
\colon I \rightarrow T^*M$ along $\gamma$ and $l\in \mathbb{N}$
satisfying $N_l=N_{l+1}$ such that for every $i=0,\dots,l$ there
exists a family of vector fields
\begin{equation}\label{Zvf} \mathscr{Z}_i(\gamma(t),u)=\{Z_1(\gamma(t),u),\dots,
Z_{n_i}(\gamma(t),u)\}\subseteq
\mathscr{V}^{i+1}_{\Sigma_2}(\gamma(t))\end{equation} fulfilling
that
\begin{enumerate} \item the sequence of constraint submanifolds $N_i$ obtained from the
presymplectic constraint algorithm is given by
\begin{equation}\label{Eq-Prop-Ni-affine}N_i=\{(\beta,u) \in N_{i-1}\mid \langle
\beta, Z_j(x,u)\rangle =0, \, j=1,\dots,n_i\}, \quad \beta\in
T^*_xM, \end{equation} for $i=1,\ldots,l$;
\item and $(\lambda,u)\in N_l$.
\end{enumerate}
\end{theorem}

\proof Since $(\gamma,u)$ is optimal, Theorem \ref{weakHPMP} can be
applied and the presymplectic constraint algorithm can be used. For
$i=0$, $N_0=\{(\beta,u) \mid \langle \beta, X_c(x)\rangle=0,
\;c=1,\dots,k\}$.
 For control-affine systems the assumption of abnormality assures that Remark \ref{remark} is true provided that
  $0$ is in the
interior of the control set, i.e.
\[{\rm cone}({\rm
conv}(\sum_{c=1}^ku^cX_c(\gamma(t))\mid u\in U))\subseteq
\mathscr{V}^1_{\Sigma_2}(\gamma(t)).\] Thus the primary constraint
submanifold is defined as in (\ref{Eq-Prop-Ni-affine}). Hence
\begin{equation*} \mathscr{Z}_0(\gamma(t),u)=\{X_1(\gamma(t)), \dots, X_k(\gamma(t))\}.\end{equation*}

Inductively, assume the result holds for $i-1$ and prove it for $i$.
If $(\lambda,u)\in N_{i-1}$, then
\[\langle \lambda, Z_j(x,u)\rangle=0, \quad j=1,\dots,n_{i-1},\]
where $Z_j(\gamma(t),u)\in \mathscr{Z}_{i-1}(\gamma(t),u) \subseteq
\mathscr{V}^i_{\Sigma_2}(\gamma(t))$. In order to obtain the
following constraint submanifold we compute
\[\frac{\rm d}{{\rm d}t} \langle \lambda, Z_j(x,u)\rangle =\langle \lambda, [\xi_0,Z_j](x,u)\rangle=0.\]
In \cite[Section 4]{1977Krener} it is proved that $[\xi_0,Z_j]$ is a
variation of order $i+1$. Thus it is true that the abnormal lift
lies in a constraint submanifold $N_i$ defined as in
(\ref{Eq-Prop-Ni-affine}). \qed

Observe that for driftless control systems this result is proved
straightforward because the control system is control-linear. Then,
\[\ds{[\xi_0,Z_j]=j^2(\nu_{\xi,\tau_2})(0)=j^2\left(\Phi^{-\xi_0}_{s} \Phi^{-Z_j}_{s}
\Phi^{\xi_0}_{s} \Phi^{Z_j}_{s} \right)(\gamma(t))(0)},\] where the
juxtaposition of maps denotes composition. Thus the vector field
$[\xi_0,Z_j]$ is obtained as a high order
elementary perturbation vector at $\gamma(t)$ of order 2. 

As in PMP, the above proposition does not guarantee the uniqueness
of the momentum, there might be more than one lift to the cotangent
bundle of the optimal curve satisfying the necessary conditions.

If the controls do not appear explicitly in the constraints defining
the submanifolds in (\ref{Eq-Prop-Ni-affine}), then the supporting
hyperplane contains subspaces generated by the family of vectors
(\ref{Zvf}) in the high order tangent perturbation cone. These
conditions impose restrictions on the momentum. If the controls
appear in the constraints, then we might be able to determine the
vector field whose integral curves are the biextremals of the OCP.

Remember that to solve uniquely Hamilton's equations ``only" the
controls and the initial condition for the momentum must be
determined. After applying the presymplectic constraint algorithm it
could happen that the initial condition for the momentum is not
completely determined. If we prefer to restrict more the candidates
to be optimal, then it must be imposed that the momentum defines a
supporting hyperplane to all the infinitesimal variations and not
only to the ones appearing in (\ref{Eq-Prop-Ni-affine}).

\begin{remark} The constructions in Definitions
\ref{Def-Krener-variations} and \ref{Def-High-Pert-Vect} have been
successfully related with the presymplectic constraint algorithm.
However, for each control system, only some infinitesimal variations
appear in the constraint algorithm as the inclusion in (\ref{Zvf})
indicates. More details are given in the following section.
\end{remark}

\subsection{Example: Affine connection control systems}\label{ExHOMP}

The affine connection control systems are a class of control
mechanical systems widely-studied in the literature
\cite{2005BulloAndrewBook}. If $Q$ is the configuration manifold and
$\nabla$ is an affine connection on $Q$, the control system on $TQ$
is the following control-affine system
\[\ds{\dot{\gamma}(t)=(Z+\sum_{c=1}^k u^cY_c^V)(\gamma(t))},\]
where $\gamma\colon I\rightarrow TQ$, $Z$ is the geodesic spray
associated with $\nabla$ and $Y^V_c$ is the vertical lift of $Y_c$.

Let $\xi_0$ be the reference trajectory with control $u_0$. The
Campbell\textendash{}Baker\textendash{}Hausdorff formula \cite{CBH}
is useful to describe specifically the sets $\mathscr{Z}_j$ whose
existence has been proved in Theorem \ref{Prop-Relation-affine}.
\begin{multline*}
\ds{\mathscr{Z}_0(\gamma(t),u)=}\\\ds{\left\{\left\{j^1\left((\Phi^{\xi_i}_s
\Phi^{\xi_0}_{-s})(\gamma(t))\right)(0)\; | \;
\xi_i=Z+\sum_{k=1}^cu_0^cY_c^V+Y_i^V
\right\}_{i=1,\dots,k}\right\}}\\
\ds{=\{Y_1^V(\gamma(t)),\dots,Y_k^V(\gamma(t)) \}}.
\end{multline*} Moreover, it can be also proved that
\begin{multline*}
\ds{\mathscr{Z}_1(\gamma(t),u)=}\\\ds{\left\{\left\{j^2\left((\Phi^{\xi_0}_{-s}
\Phi^{\xi_{-i}}_{s} \Phi^{\xi_i}_s
\Phi^{\xi_0}_{-s})(\gamma(t))\right)(0)\; | \;
\xi_{\pm i}=Z+\sum_{k=1}^cu_0^cY_c^V\pm Y_i^V \right\}_{i=1,\dots,k}\right\}}\\
 \ds{=\{[\xi_0,Y_1^V](\gamma(t)),\dots,[\xi_0,Y_k^V](\gamma(t))
\}.}
\end{multline*}

Note that if $(\lambda,u)\in N_0$, then $\lambda$ annihilates
$\mathscr{Z}_0$ and
\[\langle \lambda, [\xi_0,Y_i^V]\rangle= \langle \lambda,
[Z,Y_i^V]-\sum_{c=1}^kY_i^V(u_0^c)Y_c^V\rangle= \langle \lambda,
[Z,Y_i^V]\rangle=0. \] Thus the controls are not determined in these
first two primary constraint submanifolds. Only the momenta has been
restricted because the high order tangent perturbation cone contains
the subspaces generated by $\mathscr{Z}_0(\gamma(t),u)$ and
$\mathscr{Z}_0(\gamma(t),u)$.

Note that
\begin{align*}
\ds{j^1\left((\Phi^{\xi_{-i}}_s
\Phi^{\xi_0}_{-s})(\gamma(t))\right)(0)}&\ds{=-Y_i^V(\gamma(t)),} \\
\ds{\xi_{-i}}&=\ds{Z+\sum_{k=1}^cu_0^cY_c^V-Y_i^V},\\
\ds{j^2\left((\Phi^{\xi_0}_{-s} \Phi^{\xi_{i}}_{s} \Phi^{\xi_{-i}}_s
\Phi^{\xi_0}_{-s})(\gamma(t))\right)(0)}&\ds{=-[\xi_0,Y_i^V](\gamma(t)),}
\\ \ds{\xi_{\pm i}}&=\ds{Z+\sum_{k=1}^cu_0^cY_c^V\pm Y_i^V}.
\end{align*}

\section{Conclusions and future work}

 To sum up, HOMP provides extra conditions for optimality
related with the existence of a supporting hyperplane.
Computationally speaking, it is difficult to construct explicitly
all the infinitesimal variations. In the same way, it was difficult
to deal with the condition of maximization of the Hamiltonian over
the controls and that is why it is usually replaced by the necessary
condition $\partial H /
\partial u=0$. From this weak first order necessary condition, the
presymplectic constraint algorithm gives an idea about how to
replace the existence of the supporting hyperplane in HOMP by the
conditions in Theorem \ref{Prop-Relation-affine}. These last ones
provide less information, but they establish a connection between
the weak high order maximum principle and the application of the
presymplectic constraint algorithm to optimal control theory. The
conditions obtained by this algorithm are more suitable to work with
than to try to write down all the $(\xi,\tau_2)$-infinitesimal
variations as pointed out in the following section.

This paper sets up the foundations to characterize abnormality by
means of geometric high order necessary conditions for optimality.
We have weakened Krener's high order maximum principle to Theorem
\ref{weakHPMP}, which in turn has been weakened to Theorem
\ref{Prop-Relation-affine} to obtain easier necessary conditions for
optimality to compute. Next step is to describe when these weaker
necessary conditions can also satisfy the above-mentioned stronger
necessary conditions for optimality. In particular, how it can be
detected by means of the presymplectic algorithm if the momentum,
apart from defining a supporting hyperplane to some
$(\xi,\tau_2)$-infinitesimal variations, also determines the
supporting hyperplane in Theorems \ref{HighPMP} and \ref{weakHPMP}.
As mentioned earlier, this supporting hyperplane must be indeed a
separating hyperplane of the cone and the decreasing direction of
the cost function in order not to contradict the optimality
condition. Thus, it is also necessary to study how to detect when
the supporting hyperplane is the separating hyperplane used in the
proof of Krener's Maximum Principle.

In this regard, it could be useful to be able to define some
vector-valued quadratic forms generalizing the ones considered in
\cite{2005BulloAndrewBook} to study second order conditions for
controllability. It should be analyzed how those quadratic forms
take part in the presymplectic constraint algorithm.

In \cite{MTNSCesar} some notions of controllability for
control-affine systems have been characterized for first time in the
literature in a coordinate-free way. Thus a future research line
consists of describing necessary conditions for optimality using the
notion of affine bundles and the
jet bundle theory.  %

\section*{Acknowledgements}
We acknowledge the financial support of \emph{Ministerio de
Educaci\'on y Ciencia}, Project MTM2008-00689/MTM and the Network
Project MTM2008-03606-E/. The first author acknowledges the
financial support of Comissionat per a Universitats i Recerca del
Departament d'Innovaci\'o, Universitats i Empresa of Generalitat de
Catalunya in the preparation of this paper.



\end{document}